\documentclass{amsart}
\usepackage{latexsym}
\usepackage{amsfonts}
\usepackage{amssymb}

\newtheorem{theorem}{Theorem}[section]
\newtheorem{lemma}[theorem]{Lemma}
\numberwithin{equation}{section}

\newcommand{\rank}{\mathrm{rank}}
\newcommand{\diam}{\mathrm{diam}}
\newcommand{\Lip}{\mathrm{Lip}}
\newcommand{\supp}{\mathrm{supp}}
\newcommand{\im}{\mathrm{Im}}
\newcommand{\Diff}{\mathrm{Diff}^1(M)}

\begin{document}

\title{H\"{o}lder stability of diffeomorphisms}

\author{Jinpeng An}
\address{Department of Mathematics, ETH Zurich, 8092 Zurich,
Switzerland}
\address{Current address: Department of Pure Mathematics, University of
Waterloo Waterloo, Ontario N2L 3G1, Canada}
\email{j11an@math.uwaterloo.ca}

\keywords{Axiom A, strong transversality condition, structural
stability, H\"{o}lder regularity.}

\subjclass[2000]{37C75; 37D20.}

\begin{abstract}
We prove that a $C^2$ diffeomorphism $f$ of a compact manifold $M$
satisfies Axiom A and the strong transversality condition if and
only if it is H\"{o}lder stable, that is, any $C^1$ diffeomorphism
$g$ of $M$ sufficiently $C^1$ close to $f$ is conjugate to $f$ by a
homeomorphism which is H\"{o}lder on the whole manifold.
\end{abstract}

\maketitle


\section{Introduction}

Let $M$ be a compact $C^\infty$ manifold, $\Diff$ be the group of
$C^1$ diffeomorphisms of $M$. $f\in\Diff$ is \emph{structurally
stable} if for any $g\in\Diff$ sufficiently $C^1$ close to $f$,
there is a homeomorphism $h$ of $M$ such that $g=hfh^{-1}$. Recall
that $f$ satisfies \emph{Axiom A} if the nonwandering set $\Omega$
of $f$ is hyperbolic and the set of periodic points of $f$ is dense
in $\Omega$, $f$ satisfies the \emph{strong transversality
condition} if for any two points $x,y\in\Omega$ the stable manifold
$W^s(x)$ intersects the unstable manifold $W^u(y)$ transversally. By
the Structural Stability Theorem of Robbin, Robinson, Liao and
Ma\~{n}\'{e} \cite{R1,R2,Li,Ma}, $f\in\Diff$ is structurally stable
if and only if $f$ satisfies Axiom A and the strong transversality
condition. It is also known that in this case the conjugacy $h$ can
be chosen to be H\"{o}lder on the nonwandering set $\Omega$ of $f$
(see \cite[Theorem 19.1.2]{KH}).

In this paper, we prove that in the above case, the conjugacy $h$
can be chosen to be H\"{o}lder not only on $\Omega$ but also on the
whole manifold $M$. We say that a deffeomorphism $f$ of $M$ is
\emph{H\"{o}lder stable} if for any $g\in\Diff$ sufficiently $C^1$
close to $f$, there is a H\"{o}lder homeomorphism $h$ of $M$ such
that $g=hfh^{-1}$ (This notion should not be confused with the
notion of $C^r$ structural stability of a $C^r$ diffeomorphism, for
which $g$ is $C^r$ close to $f$ and the conjugacy $h$ is only
required to be continuous). We prove that Axiom A plus the strong
transversality condition is also equivalent to H\"{o}lder stability.
For simplicity, we assume that $f$ is $C^2$.

\begin{theorem}\label{T:main}
Let $f$ be a $C^2$ diffeomorphism of a compact $C^\infty$ manifold
$M$. Then $f$ is H\"{o}lder stable if and only if $f$ satisfies
Axiom A and the strong transversality condition.
\end{theorem}

Since H\"{o}lder stability implies structural stability, to prove
Theorem \ref{T:main}, it is sufficient by the Structural Stability
Theorem to prove that Axiom A plus the strong transversality
condition implies H\"{o}lder stability.

To state the quantitative result, we recall the notion of
hyperbolicity. The nonwandering set $\Omega$ of a diffeomorphism $f$
is \emph{hyperbolic} if the restriction $TM|_{\Omega}$ of the
tangent bundle $TM$ on $\Omega$ admits a $Tf$-invariant continuous
splitting $TM|_{\Omega}=E^u\oplus E^s$ such that for some
$\lambda\in(0,1)$,
\begin{equation}\label{E:lambda}
\|Tf^{-1}|_{E^u}\|\leq \lambda, \quad \|Tf|_{E^s}\|\leq \lambda.
\end{equation}
Here the norm is evaluated with respect to some adapted smooth
Riemannian metric on $M$.

\begin{theorem}\label{T:main2}
Let $f$ be a $C^2$ diffeomorphism of a compact $C^\infty$ manifold
$M$ satisfying Axiom A and the strong transversality condition. Let
$\lambda\in(0,1)$ be as in \eqref{E:lambda},
$l=\max\{\Lip(f),\Lip(f^{-1})\}$. Suppose $\alpha\in(0,1)$ satisfies
$\lambda l^\alpha<1$. Then for any $C^\alpha$ neighborhood
$\mathcal{V}$ of the identity map in $C^\alpha(M,M)$, there exists a
$C^1$ neighborhood $\mathcal{N}$ of $f$ in $\Diff$ such that for
every $g\in\mathcal{N}$, there is a homeomorphism $h$ of $M$ in
$\mathcal{V}$ such that $g=hfh^{-1}$, and the assignment $g\mapsto
h$ is $C^1$ as a map $\mathcal{N}\rightarrow C^0(M,M)$ and sends $f$
to the identity.
\end{theorem}

Here $\Lip(f)$ denotes the Lipschitz constant of $f$,
$C^\alpha(M,M)$ and $C^0(M,M)$ are the Banach manifolds of
$C^\alpha$ and $C^0$ maps on $M$, respectively.

H\"{o}lder stabillity over hyperbolic sets is well known
(\cite[Theorem 19.1.2]{KH}). It is also well known that the
(un)stable distributions and (un)stable foliations over hyperbolic
sets are H\"{o}lder continuous (\cite[Section 19.1]{KH}). For more
results on H\"{o}lder regularity for hyperbolic dynamical systems,
see \cite[Section 2.3]{Ha}.

One can not expect more regularity of the conjugacy $h$ than to be
H\"{o}lder. For example, Lipschitz conjugacies almost never exist.
But for dynamical systems of large group actions, $C^r$ or
$C^\infty$ conjugacies may exist (see \cite{Fi} and the references
therein).

Our proof of Theorem \ref{T:main2} follows the approach of
Robbin-Robinson \cite{R1,R2}, where the result that Axiom A plus the
strong transversality condition implies structural stability is
proved. As in Robbin \cite{R1}, we divide the proof into three
steps, which are the contents of the following three sections.

In Section 2, we prove that for each component $\Omega_i$ in the
spectral decomposition of $\Omega$, the splitting
$TM|_{\Omega_i}=E^u|_{\Omega_i}\oplus E^s|_{\Omega_i}$ can be
extended to a $Tf$-invariant splitting
$TM|_{\mathcal{O}(U_i)}=E^u_i\oplus E^s_i$ satisfying certain
compatibility condition, where $U_i$ is a neighborhood of
$\Omega_i$,
$\mathcal{O}(U_i)=\bigcup_{n=-\infty}^{+\infty}f^n(U_i)$. The proof
follows ideas in \cite{R1,R2}. But since we require that the
extended splitting to be H\"{o}lder, and the metric $d$ on $M$,
unlike Robbin's metric $d_f$ \cite{R1}, is not $f$-preserving, we
need more careful topological arguments. Indeed, we can only prove
that the extended bundles $E^u_i$ and $E^s_i$ are H\"{o}lder on
$\bigcup_{n=-N}^{N}f^n(U_i)$ for every $N>0$. But this is sufficient
for us to derive further results. In Section 2 we only need the
weaker restriction $\lambda^2l^\alpha<1$ on the H\"{o}lder exponent
$\alpha$ comparing with Theorem \ref{T:main2}, and the case of
$\alpha=1$ is allowed, which means as usual that the subbundles are
Lipschitz.

Using the extended splitting in Section 2, we prove in Section 3
that the induced operator $f_\sharp$ of $f$ on the Banach space of
$C^0$ vector fields has a right inverse which restricts to a
continuous linear operator on the Banach space of $C^\alpha$ and
$d_f$-Lipschitz vector fields. The proof is also motivated by
\cite{R1}. But as in Section 2, since $f$ does not preserve the
metric $d$, some different topological arguments are needed. The
condition of $\alpha\neq1$ is not explicitly used in the proof. But
since it is easy to see that $l\geq\lambda^{-1}$, the inequality
$\lambda l^\alpha<1$ for $\alpha=1$ never holds. So the case of
$\alpha=1$ is automatically excluded.

In Section 4 we finish the proof of Theorem \ref{T:main2}. We first
prove a version of Implicit Function Theorem for Banach spaces
involving non-closed subspaces. Then using the result in Section 3,
we can apply the Implicit Function Theorem to the $C^1$ map
$\Psi:\Diff\times C^0(M,M)\rightarrow C^0(M,M)$,
$\Psi(g,h)=ghf^{-1}$ to obtain a fixed point $h$ of $\Psi(g,\cdot)$
for $g$ sufficiently $C^1$ close to $f$, and $h$ is sufficiently
$C^\alpha$ and $d_f$-Lipschitz close to the identity. As in
\cite{R1,R2}, the fact that $h$ is $d_f$-Lipschitz close to the
identity implies $h$ is a homeomorphism.

Most arguments concerning $C^0$ estimates in this paper are borrowed
from \cite{R1,R2} except for a few changes of details. But to
introduce notations in order to perform the $C^\alpha$ estimates, it
seems necessary to repeat some of them.

The author would like to thank Professors Boris Hasselblatt and Lan
Wen for useful comments.

\section{Extensions of the splitting}

In this section we prove that the splitting $TM|_{\Omega}=E^u\oplus
E^s$ can be extended to a neighborhood of each component of $\Omega$
and satisfies certain compatibility condition. This is motivated by
\cite[Theorem 8.4, C]{R1} and \cite[Theorem 3.1, 5.1]{R2}.

We first collect some standard facts that are used in the proof of
Theorem \ref{T:1} below. Most of them can be found in
\cite{HPPS,R1,S}. Let $f$ be a diffeomorphism of a compact manifold
$M$ satisfying Axiom A and the strong transversality condition. Let
$\Omega=\Omega_1\cup\cdots\cup\Omega_k$ be the spectral
decomposition of the nonwandering set $\Omega$ of $f$. Each
$\Omega_i$ is a closed topological transitive hyperbolic
$f$-invariant subset of $M$, and $E^u, E^s$ have constant ranks on
$\Omega_i$. The components $\Omega_i$ can be ordered in such a way
that $i<j$ implies $W^s(\Omega_i)\cap W^u(\Omega_j)=\emptyset$,
where $W^\sigma(\Omega_i)=\bigcup_{x\in\Omega_i}W^\sigma(x)$,
$\sigma=u,s$. For a subset $U$ of $M$, denote
$\mathcal{O}(U)=\bigcup_{n=-\infty}^{+\infty}f^n(U)$,
$\mathcal{O}^+(U)=\bigcup_{n=0}^{+\infty}f^n(U)$,
$\mathcal{O}^-(U)=\bigcup_{n=0}^{+\infty}f^{-n}(U)$. Then for
$\Omega_i, \Omega_j$ such that $W^s(\Omega_i)\cap
W^u(\Omega_j)=\emptyset$ and sufficiently small neighborhoods $U_i,
U_j$ of $\Omega_i$ and $\Omega_j$,
$\mathcal{O}^-(U_i)\cap\mathcal{O}^+(U_j)=\emptyset$. A subset $U$
of $M$ is called \emph{unrevisited} if $x\in U$, $n>0$, $f^n(x)\in
U$ imply $f^m(x)\in U$ for $0<m<n$. Then each $\Omega_i$ has
arbitrarily small unrevisited open neighborhood.

We fix an $\Omega_i$. For $x\in\Omega_i$ and $\delta>0$, let
$W^u_\delta(x)$ and $W^s_\delta(x)$ be the local unstable and stable
manifolds of size $\delta$ at $x$. Let
$W^\sigma_\delta(\Omega_i)=\bigcup_{x\in\Omega_i}W^\sigma_\delta(x)$,
$\sigma=u,s$. For $\delta$ sufficiently small,
$W^\sigma_\delta(\Omega_i)$ has arbitrarily small unrevisited open
neighborhood. Let $D=\overline{W^s_\delta(\Omega_i)\setminus
f(W^s_\delta(\Omega_i))}$. Then for $\delta$ sufficiently small, $D$
has arbitrarily small unrevisited open neighborhood, and for any
open neighborhood $Q$ of $D$, the set
$W^u(\Omega_i)\cup\mathcal{O}^+(Q)$ is an unrevisited open
neighborhood of $\Omega_i$.

As in \cite{R1,R2}, we introduce the metric $d_f$ on $M$ by
$d_f(x,y)=\sup_{n\in\mathbb{Z}}d(f^n(x),d^n(y))$, where $d$ is the
metric induced from some Riemannian metric on $M$.

\begin{theorem}\label{T:1}
Let $f$ be a $C^2$ diffeomorphism of $M$ satisfying Axiom A and the
strong transversality condition,
$\Omega=\Omega_1\cup\cdots\cup\Omega_k$ be the spectral
decomposition, and the components $\Omega_i$ are ordered as above.
Let $\lambda\in(0,1)$ be as in \eqref{E:lambda},
$l=\max\{\Lip(f),\Lip(f^{-1})\}$. Suppose $\alpha\in(0,1]$ satisfies
$\lambda^2 l^\alpha<1$. Then for any $\lambda'\in(\lambda,1)$, there
exist for each $1\leq i\leq k$ an open neighborhood $U_i$ of
$\Omega_i$ and two $Tf$-invariant continuous subbundles $E^\sigma_i$
of
$TM|_{\mathcal{O}(U_i)}$, $\sigma=u,s$, such that\\
(i) $TM|_{\mathcal{O}(U_i)}=E^u_i\oplus E^s_i$;\\
(ii) $E^\sigma_i$ is $C^\alpha$ and $d_f$-Lipschitz on
$\bigcup_{n=-N}^{N}f^n(U_i)$ for every $N>0$;\\
(iii) $\|Tf^{-1}|_{(E^u_i)_x}\|\leq\lambda'$,
$\|Tf|_{(E^s_i)_x}\|\leq\lambda'$
for  $x\in U_i$;\\
(iv) for $i<j$,
$\mathcal{O}^-(U_i)\cap\mathcal{O}^+(U_j)=\emptyset$, and
$(E^s_i)_x\subset(E^s_j)_x$, $(E^u_j)_x\subset(E^u_i)_x$ for every
$x\in\mathcal{O}^+(U_i)\cap\mathcal{O}^-(U_j)$.
\end{theorem}

\begin{proof}
We extend the definition of the bundles $E^u$ and $E^s$ on $\Omega$
as $E^u=\{v\in TM:\lim_{n\rightarrow+\infty}|Tf^{-n}(v)|=0\}$,
$E^s=\{v\in TM:\lim_{n\rightarrow+\infty}|Tf^n(v)|=0\}$, and denote
$E^\sigma_x=E^\sigma\cap T_xM$, $\sigma=u,s$. By the strong
transversality condition, $T_xM=E^u_x+E^s_x$ for every $x\in M$. For
each $\Omega_i$, $E^\sigma|_{W^\sigma(\Omega_i)}$ is a continuous
subbundle of $TM|_{W^\sigma(\Omega_i)}$ with constant rank.

As in \cite[Section 10]{R1}, to prove Theorem \ref{T:1}, it is
sufficient to prove that under the conditions of Theorem \ref{T:1},
there exist for each $i$ an open neighborhood $U_i$ of $\Omega_i$
and a $Tf$-invariant continuous subbundle $E^u_i$ of
$TM|_{\mathcal{O}(U_i)}$
such that\\
(i') $E^u_i|_{\Omega_i}=E^u|_{\Omega_i}$;\\
(ii') $E^u_i$ is $C^\alpha$ and $d_f$-Lipschitz on
$\bigcup_{n=-N}^{N}f^n(U_i)$ for every $N>0$;\\
(iii') for $i<j$,
$\mathcal{O}^-(U_i)\cap\mathcal{O}^+(U_j)=\emptyset$, and
$(E^u_j)_x\subset(E^u_i)_x$ for every
$x\in\mathcal{O}^+(U_i)\cap\mathcal{O}^-(U_j)$;\\
(iv') $T_xM=(E^u_i)_x+E^s_x$ for every $x\in\mathcal{O}(U_i)$.

We prove this by induction on $i=1,\cdots,k$. Let $1\leq i\leq k$.
Suppose that for $j<i$, $U_j$ and $E^u_j$ have been defined and
satisfy (i')--(iv') (for $i=1$ nothing is defined). We construct
$U_i$ and $E^u_i$ satisfying (i')--(iv').

Let $\lambda<\lambda_1<\lambda_2<\lambda_3<1$ be such that
$\lambda_3^2l^\alpha<1$. Let $V_1$ be an open neighborhood of
$\Omega_i$ such that
$\mathcal{O}^+(V_1)\cap\mathcal{O}^-(U_j)=\emptyset$ for all $j<i$
(shrinking $U_j$, $j<i$ if necessary). Choose continuous subbundles
$\widetilde{E}^u, \widetilde{E}^s$ of $TM|_{V_1}$ with
$\widetilde{E}^\sigma|_{V_1\cap
W^\sigma(\Omega_i)}=E^\sigma|_{V_1\cap W^\sigma(\Omega_i)}$,
$\sigma=u,s$. Since $TM|_{\Omega_i}=E^u|_{\Omega_i}\oplus
E^s|_{\Omega_i}$, shrinking $V_1$ if necessary, we may assume that
$TM|_{V_1}=\widetilde{E}^u\oplus \widetilde{E}^s$. Write
$Tf|_{V_1\cap f^{-1}(V_1)}$ as
$$T_xf=\begin{pmatrix}
\widetilde{F}^{uu}_x & \widetilde{F}^{su}_x\\
\widetilde{F}^{us}_x & \widetilde{F}^{ss}_x
\end{pmatrix}$$ with respect to the splitting $TM|_{V_1}=\widetilde{E}^u\oplus
\widetilde{E}^s$, $x\in V_1\cap f^{-1}(V_1)$. Since
$\|(\widetilde{F}^{uu}_x)^{-1}\|\leq\lambda$,
$\|\widetilde{F}^{ss}_x\|\leq\lambda$, $\widetilde{F}^{us}_x=0$ for
$x\in\Omega_i$, by making $V_1$ smaller, we may assume that
$\|(\widetilde{F}^{uu}_x)^{-1}\|\leq\lambda_1$,
$\|\widetilde{F}^{ss}_x\|+\|\widetilde{F}^{us}_x\|\leq\lambda_1$ for
$x\in V_1\cap f^{-1}(V_1)$. Note that since
$\widetilde{E}^s|_{V_1\cap W^s(\Omega_i)}=E^s|_{V_1\cap
W^s(\Omega_i)}$ and $Tf(E^s)=E^s$, $\widetilde{F}^{su}_x|_{V_1\cap
W^s(\Omega_i)}=0$.

Choose $\delta>0$ such that $W^s_{\delta}(\Omega_i)\subset V_1$, and
such that $W^s_{\delta}(\Omega_i)$ has arbitrarily small unrevisited
open neighborhood. Let $D=\overline{W^s_\delta(\Omega_i)\setminus
f(W^s_\delta(\Omega_i))}$. Similar to the arguments in \cite[page
488--491]{R1}, we can prove (after possibly shrinking of $U_j, j<i$
in the induction hypothesis) that there exist an open neighborhood
$Q_1\subset V_1$ of $D$ and a $C^\alpha$ and $d_f$-Lipschitz
subbundle $E^u_{i0}$ of $TM|_{Q_1}$ such that\\
(1) $Tf((E^u_{i0})_x)=(E^u_{i0})_{f(x)}$ for $x\in Q_1\cap
f^{-1}(Q_1)$;\\
(2) $T_xM=(E^u_{i0})_x\oplus \widetilde{E}^s_x$ and
$T_xM=(E^u_{i0})_x+E^s_x$ for $x\in Q_1$;\\
(3) $(E^u_{i0})_x\subset (E^u_j)_x$ if $j<i$ and $x\in Q_1\cap
\mathcal{O}^+(U_j)$.\\
We may also assume that $Q_1\cap f^2(Q_1)=\emptyset$.

Since $TM|_{Q_1}=E^u_{i0}\oplus\widetilde{E}^s|_{Q_1}$, there exists
a continuous vector bundle morphism
$\widetilde{\tau}_0:\widetilde{E}^u|_{Q_1}\rightarrow\widetilde{E}^s|_{Q_1}$
such that $E^u_{i0}=\im(id,\widetilde{\tau}_0)$. By making $Q_1$
smaller, we may assume that $\|\widetilde{\tau}_0\|$ is bounded, say
$\|\widetilde{\tau}_0\|\leq\frac{r}{2}$ for some $r\geq1$.

Choose $\varepsilon>0$ such that $r\varepsilon
\leq\lambda_2^{-1}-\lambda_3^{-1}$. Since
$\widetilde{F}^{su}_x|_{W^s_{\delta}(\Omega_i)}=0$, we may choose an
unrevisited open neighborhood $V_2\subset V_1$ of
$W^s_{\delta}(\Omega_i)$ such that
$\|\widetilde{F}^{su}_x\|\leq\frac{\varepsilon}{2}$ for $x\in V_2$.
Let $Q_2\subset Q_1\cap V_2$ be an unrevisited open neighborhood of
$D$. Choose $C^1$ approximations $\bar{E}^u, \bar{E}^s$ of
$\widetilde{E}^u|_{V_2}, \widetilde{E}^s|_{V_2}$ such that\\
(1) $TM|_{V_2}=\bar{E}^u\oplus\bar{E}^s$, and if
$$T_xf=\begin{pmatrix}
F^{uu}_x & F^{su}_x\\
F^{us}_x & F^{ss}_x
\end{pmatrix}$$ with respect to this splitting, then
$\|F^{su}_x\|\leq\varepsilon$, $\|(F^{uu}_x)^{-1}\|\leq\lambda_2$,
$\|F^{ss}_x\|+\|F^{us}_x\|\leq\lambda_2$ for $x\in
V_2\cap f^{-1}(V_2)$;\\
(2) $TM|_{Q_2}=E^u_{i0}|_{Q_2}\oplus\bar{E}^s|_{Q_2}$, and if
$\tau_0:\bar{E}^u|_{Q_2}\rightarrow\bar{E}^s|_{Q_2}$ is the vector
bundle morphism
such that $E^u_{i0}|_{Q_2}=\im(id,\tau_0)$, then $\|\tau_0\|\leq r$;\\
(3) there exists a continuous vector bundle morphism
$\tau'_0:\bar{E}^u|_{V_2\cap
W^u(\Omega_i)}\rightarrow\bar{E}^s|_{V_2\cap W^u(\Omega_i)}$ such
that $E^u|_{V_2\cap
W^u(\Omega_i)}=\im(id,\tau'_0)$, and $\|\tau'_0\|\leq r$.\\
Note that since $f$ is $C^2$ and the splitting
$TM|_{V_2}=\bar{E}^u\oplus\bar{E}^s$ is $C^1$, $F^{\sigma\sigma'}$
is $C^1$, where $\sigma,\sigma'=u,s$. Note also that since
$E^u_{i0}$ is $C^\alpha$ and $d_f$-Lipschitz, $\tau_0$ is $C^\alpha$
and $d_f$-Lipschitz.

Consider the smooth vector bundle $\mathcal{L}$ over $V_2$ whose
fiber $\mathcal{L}_x$ at $x\in V_2$ is the space
$\mathcal{L}(\bar{E}^u_x, \bar{E}^s_x)$ of linear maps from
$\bar{E}^u_x$ to $\bar{E}^s_x$. A section $\tau$ of $\mathcal{L}$ is
a vector bundle morphism from $\bar{E}^u$ to $\bar{E}^s$ covering
the identity. Let $\mathcal{L}(r)_x$ be the disc
$\{g\in\mathcal{L}_x:\|g\|\leq r\}$ in $\mathcal{L}_x$, and
$\mathcal{L}(r)=\bigcup_{x\in V_2}\mathcal{L}(r)_x$ be the disc
bundle of $\mathcal{L}$. Let $x\in V_2\cap f^{-1}(V_2)$. For
$g\in\mathcal{L}(r)_x$, define
$$\varphi^\sigma_g=F^{u\sigma}_x+F^{s\sigma}_xg\in\mathcal{L}
(\bar{E}^u_x, \bar{E}^\sigma_{f(x)}),$$ $\sigma=u,s$.
Then for $v\in\bar{E}^u_x$, we have
$$|\varphi^u_g(v)|=|F^{uu}(v)+F^{su}g(v)|\geq|F^{uu}(v)|-|F^{su}g(v)|\geq(\lambda_2^{-1}-r\varepsilon)|v|.$$
So $\varphi^u_g$ is invertible, and
$$\|(\varphi^u_g)^{-1}\|\leq(\lambda_2^{-1}-r\varepsilon)^{-1}\leq\lambda_3.$$ we
also have
$$\|\varphi^s_g\|=\|F^{us}\|+\|F^{ss}g\|\leq\|F^{us}\|+r\|F^{ss}\|\leq\lambda_2r.$$
Thus if we define the graph transform of $g\in\mathcal{L}(r)_x$ by
$$\Gamma(g)=\varphi^s_g(\varphi^u_g)^{-1}\in\mathcal{L}_{f(x)},$$ then
\begin{equation}\label{E:<=K}
\|\Gamma(g)\|\leq\lambda_3\lambda_2r\leq r.
\end{equation}
Note that since $F^{\sigma\sigma'}$ is $C^1$, the map
$\Gamma:\mathcal{L}(r)|_{V_2\cap
f^{-1}(V_2)}\rightarrow\mathcal{L}(r)|_{f(V_2)\cap V_2}$ is $C^1$.

Let $x\in V_2\cap f^{-1}(V_2)$, $g_1, g_2\in\mathcal{L}(r)_x$. Then
$$\|\varphi^u_{g_1}-\varphi^u_{g_2}\|=\|F^{su}(g_1-g_2)\|\leq\varepsilon\|(g_1-g_2)\|,$$
$$\|\varphi^s_{g_1}-\varphi^s_{g_2}\|=\|F^{ss}(g_1-g_2)\|\leq\lambda_2\|(g_1-g_2)\|.$$
Hence
\begin{align}\label{E:contraction}
&\notag\|\Gamma(g_1)-\Gamma(g_2)\|\\
\leq&\notag\|\varphi^s_{g_1}(\varphi^u_{g_1})^{-1}-\varphi^s_{g_1}(\varphi^u_{g_2})^{-1}\|
+\|\varphi^s_{g_1}(\varphi^u_{g_2})^{-1}-\varphi^s_{g_2}(\varphi^u_{g_2})^{-1}\|\\
\leq&\notag\|\varphi^s_{g_1}\|
\|(\varphi^u_{g_1})^{-1}\varphi^u_{g_2}(\varphi^u_{g_2})^{-1}-
(\varphi^u_{g_1})^{-1}\varphi^u_{g_1}(\varphi^u_{g_2})^{-1}\|\\
&\notag+\|(\varphi^u_{g_2})^{-1}\| \|\varphi^s_{g_1}-\varphi^s_{g_2}\|\\
\leq&\|\varphi^s_{g_1}\|\|(\varphi^u_{g_1})^{-1}\|\|(\varphi^u_{g_2})^{-1}\|
\|\varphi^u_{g_2}-\varphi^u_{g_1}\|+\|(\varphi^u_{g_2})^{-1}\|
\|\varphi^s_{g_1}-\varphi^s_{g_2}\|\\
\leq&\notag\frac{\lambda_2r}{(\lambda_2^{-1}-r\varepsilon)^2}\|\varphi^u_{g_1}-\varphi^u_{g_2}\|+
\frac{1}{\lambda_2^{-1}-r\varepsilon}\|\varphi^s_{g_1}-\varphi^s_{g_2}\|\\
\leq&\notag(\frac{\lambda_2r\varepsilon}{(\lambda_2^{-1}-r\varepsilon)^2}+
\frac{\lambda_2}{\lambda_2^{-1}-r\varepsilon})\|g_1-g_2\|\\
\leq&\notag\lambda_3^2\|g_1-g_2\|.
\end{align}

For the convenience of the following discussion, we embed $M$
isometrically into some Euclidian space $\mathbb{R}^N$. Then for
$x\in V_2$, $\bar{E}^s_x$ and $\bar{E}^u_x$ can be viewed as
subspaces of $\mathbb{R}^N$, and we have the identification
$$\mathcal{L}_x=\mathcal{L}(\bar{E}^u_x,\bar{E}^s_x)\cong
\{g\in\mathcal{L}(\mathbb{R}^N,\mathbb{R}^N)|\bar{E}^s_x\oplus
T_xM^\bot\subset\ker(g), \im(g)\subset\bar{E}^s_x\}.$$ Then for
$g_1,g_2\in \mathcal{L}$ with different base points, the summation
$g_1+g_2$ and its norm $\|g_1+g_2\|$ make sense, as they are viewed
as elements in $\mathcal{L}(\mathbb{R}^N,\mathbb{R}^N)$. Let
$\Gamma_x=\Gamma|_{\mathcal{L}_x}$ be the restriction of $\Gamma$ on
the fiber $\mathcal{L}_x$. Since the map
$\Gamma:\mathcal{L}(r)|_{V_2\cap
f^{-1}(V_2)}\rightarrow\mathcal{L}(r)|_{f(V_2)\cap V_2}$ is $C^1$,
it is Lipschitz and $C^\alpha$, which means that there exists $C>0$
such that
\begin{equation}\label{E:Lip}
\|\Gamma_x(g_1)-\Gamma_y(g_2)\|\leq C\min\{\|g_1-g_2\|+d(x,y),
(\|g_1-g_2\|+d(x,y))^\alpha\}
\end{equation}
for any $x,y\in V_2\cap f^{-1}(V_2)$ and $g_1\in\mathcal{L}(r)_x,
g_2\in\mathcal{L}(r)_y$. Note that since $\Gamma$ covers $f$ and $f$
is Lipschitz, we have indeed omitted a term $d(f(x),f(y))$ in the
left hand side of \eqref{E:Lip}.

Recall that $D\cap W^u(\Omega_i)=\emptyset$. So there exist $d_0>0$
and an unrevisited open neighborhood $Q_3\subset Q_2$ of $D$ such
that $d(Q_3,W^u(\Omega_i))\geq d_0$, and such that $x\in Q_3, y\in
M, d(x,y)<d_0$ imply $y\in Q_2$. Let $V_3=V_2\cap(W^u(\Omega_i)\cup
\mathcal{O}^+(Q_3))$. Since $V_2$ and $W^u(\Omega_i)\cup
\mathcal{O}^+(Q_3)$ are unrevisited open neighborhoods of
$\Omega_i$, so is $V_3$.

To simplify notations, we denote $$\rho_f(x,y)=\min\{d(x,y)^\alpha,
d_f(x,y)\}$$ for $x,y\in M$. Then a section $\tau$ of $\mathcal{L}$
is $C^\alpha$ and $d_f$-Lipschitz if and only if $$\sup_{x,y\in
V_2,x\neq y}\frac{\|\tau(x)-\tau(y)\|}{\rho_f(x,y)}<+\infty.$$ Now
we choose
$$K\geq\max\left\{\frac{2r\diam(M)^{1-\alpha}}{d_0}, \frac{2r}{d_0}, \frac{Cl^\alpha}
{1-\lambda_3^2l^\alpha}, \frac{C}{1-\lambda_3^2}\right\}$$ such that
$$\|\tau_0(x)-\tau_0(y)\|\leq K\rho_f(x,y)$$ for $x,y\in Q_2$, where $\diam(M)$ is the diameter of $M$.
Let
$$
\Sigma=\{\text{continuous sections $\tau$ of $\mathcal{L}(r)|_{V_3}:
\|\tau(x)-\tau(y)\|\leq K\rho_f(x,y), \tau|_{Q_3}=\tau_0|_{Q_3}$}\}.
$$
$\Sigma$ is a closed subset of the Banach space of continuous
bounded sections of $\mathcal{L}|_{V_3}$. By taking a bump function
on $M$ which is $1$ in $Q_3$ and $0$ outside $Q_2$ and enlarging $K$
if necessary, it is easy to see that $\Sigma$ is nonempty. Define
the graph transform $F_\sharp(\tau)$ of $\tau\in\Sigma$ as the
section
$$F_\sharp(\tau)(x)=\begin{cases}
\Gamma(\tau(f^{-1}(x))), &x\in f(V_3)\cap V_3;\\
\tau(x), &x\in V_3\setminus f(V_3)
\end{cases}$$ of $\mathcal{L}|_{V_3}$.
We prove that $F_\sharp$ maps $\Sigma$ into $\Sigma$ and is a
contraction on $\Sigma$.

First we show that $V_3=(f(V_3)\cap V_3)\cup Q_3$. Let $x\in V_3$.
Recall that $V_3=V_2\cap(W^u(\Omega_i)\cup \mathcal{O}^+(Q_3))$. If
$x\in V_2\cap W^u(\Omega_i)$, then there exists $n\geq1$ such that
$f^{-n}(x)\in V_2$. Since $V_2$ is unrevisited, $f^{-1}(x)\in V_2$.
We also have $f^{-1}(x)\in W^u(\Omega_i)$. Hence $x\in f(V_3)\cap
V_3$. If $x\in V_2\cap\mathcal{O}^+(Q_3)$, there exists $y\in Q_3$
such that $x=f^n(y)$ for some $n\geq0$. If $n=0$ then $x\in Q_3$. If
$n\geq1$, since $V_2$ is unrevisited, $f^{n-1}(y)\in V_2$. Hence
$x\in f(V_3)\cap V_3$. This proves $V_3=(f(V_3)\cap V_3)\cup Q_3$.

Let $\tau\in\Sigma$. We show that
$F_\sharp(\tau)|_{Q_3}=\tau_0|_{Q_3}$. Let $x\in Q_3$. If $x\in
V_3\setminus f(V_3)$, then $F_\sharp(\tau)(x)=\tau(x)=\tau_0(x)$. If
$x\in Q_3\setminus(V_3\setminus f(V_3))=Q_3\cap f(V_3)$, then
$f^{-n}(x)\in Q_3$ for some $n\geq1$. Since $Q_3$ is unrevisited,
$f^{-1}(x)\in Q_3$. So
$F_\sharp(\tau)(x)=\Gamma(\tau(f^{-1}(x)))=\Gamma(\tau_0(f^{-1}(x)))=\tau_0(x)$.
So $F_\sharp(\tau)|_{Q_3}=\tau_0|_{Q_3}$.

Now $F_\sharp(\tau)$ is continuous on $Q_3$ and $f(V_3)\cap V_3$.
Since $f(V_3)\cap V_3$ and $Q_3$ are open in $V_3$ and
$V_3=(f(V_3)\cap V_3)\cup Q_3$, $F_\sharp(\tau)$ is continuous on
$V_3$.

By \eqref{E:<=K}, $\|F_\sharp(\tau)(x)\|\leq r$ for $x\in f(V_3)\cap
V_3$. So $\|F_\sharp(\tau)\|\leq r$.

Now we show that $\|F_\sharp(\tau)(x)-F_\sharp(\tau)(y)\|\leq
K\rho_f(x,y)$ for $\tau\in\Sigma$ and $x,y\in V_3$. There are three
cases.

(1) $x,y\in V_3\setminus f(V_3)$. This is obvious since
$F_\sharp(\tau)|_{Q_3}=\tau_0|_{Q_3}$ and $V_3\setminus
f(V_3)\subset Q_3$.

(2) $x\in V_3\setminus f(V_3), y\in f(V_3)\cap V_3$. If $d(x,y)\geq
d_0$, then
\begin{align*}
&\|F_\sharp(\tau)(x)-F_\sharp(\tau)(y)\|\leq
2r\leq\frac{2r}{d_0}d(x,y)\\
\leq&\frac{2r\max\{\diam(M)^{1-\alpha},1\}}{d_0}\rho_f(x,y)\leq
K\rho_f(x,y).
\end{align*}
Suppose $d(x,y)<d_0$. Since $x\in Q_3$, we have $y\in Q_2$ and
$y\notin W^u(\Omega_i)$. So there exists $n\geq1$ such that
$f^{-n}(y)\in Q_3$. But $Q_2$ is unrevisited and $Q_2\cap
f^2(Q_2)=\emptyset$. So we must have $n=1$ and then
$F_\sharp(\tau)(y)=\Gamma(\tau(f^{-1}(y)))=\Gamma(\tau_0(f^{-1}(y)))=\tau_0(y)$.
So
$\|F_\sharp(\tau)(x)-F_\sharp(\tau)(y)\|=\|\tau_0(x)-\tau_0(y)\|\leq
K\rho_f(x,y)$.

(3) $x,y\in f(V_3)\cap V_3$. By \eqref{E:contraction} and
\eqref{E:Lip},
\begin{align*}
&\|F_\sharp(\tau)(x)-F_\sharp(\tau)(y)\|\\
=&\|\Gamma_{f^{-1}(x)}(\tau(f^{-1}(x)))-\Gamma_{f^{-1}(y)}(\tau(f^{-1}(y)))\|\\
\leq&\|\Gamma_{f^{-1}(x)}(\tau(f^{-1}(x)))-\Gamma_{f^{-1}(x)}(\tau(f^{-1}(y)))\|\\
&+\|\Gamma_{f^{-1}(x)}(\tau(f^{-1}(y)))-\Gamma_{f^{-1}(y)}(\tau(f^{-1}(y)))\|\\
\leq&\lambda_3^2\|\tau(f^{-1}(x))-\tau(f^{-1}(y))\|+Cd(f^{-1}(x),f^{-1}(y))^\alpha\\
\leq&\lambda_3^2Kd(f^{-1}(x),f^{-1}(y))^\alpha+Cd(f^{-1}(x),f^{-1}(y))^\alpha\\
\leq&(\lambda_3^2K+C)l^\alpha d(x,y)^\alpha\\
\leq&Kd(x,y)^\alpha.
\end{align*}
Similarly,\begin{align*}
&\|F_\sharp(\tau)(x)-F_\sharp(\tau)(y)\|\\
\leq&\lambda_3^2Kd_f(f^{-1}(x),f^{-1}(y))+Cd(f^{-1}(x),f^{-1}(y))\\
\leq&(\lambda_3^2K+C)d_f(x,y)\\
\leq&Kd_f(x,y).
\end{align*}
So $\|F_\sharp(\tau)(x)-F_\sharp(\tau)(y)\|\leq K\rho_f(x,y)$.

This proves that $F_\sharp$ maps $\Sigma$ into $\Sigma$. By
\eqref{E:contraction}, $F_\sharp$ is a contraction on $\Sigma$. So
there is a fixed point $\bar{\tau}$ of $F_\sharp$ in $\Sigma$.

Choose $\delta'>0$ such that $W^u_{\delta'}(\Omega_i)\subset
f(V_3)\cap V_3$. We prove that
$\bar{\tau}|_{W^u_{\delta'}(\Omega_i)}=\tau_0'|_{W^u_{\delta'}(\Omega_i)}$.
Let $x_0\in W^u_{\delta'}(\Omega_i)$ be such that
$\|\bar{\tau}|_{W^u_{\delta'}(\Omega_i)}(x)-\tau_0'|_{W^u_{\delta'}(\Omega_i)}(x)\|$
assumes maximal value at $x_0$. Then
\begin{align*}
&\|\bar{\tau}(x_0)-\tau_0'(x_0)\|\\
=&\|\Gamma(\bar{\tau}(f^{-1}(x_0)))-\Gamma(\tau_0'(f^{-1}(x_0)))\|\\
\leq&\lambda_3^2\|\bar{\tau}(f^{-1}(x_0))-\tau_0'(f^{-1}(x_0))\|\\
\leq&\lambda_3^2\|\bar{\tau}(x_0)-\tau_0'(x_0)\|.
\end{align*}
Hence $(1-\lambda_3^2)\|\bar{\tau}(x_0)-\tau_0'(x_0)\|\leq0$, which
implies that $\|\bar{\tau}(x_0)-\tau_0'(x_0)\|=0$.

Define the $C^\alpha$ and $d_f$-Lipschitz subbundle $E^u_{i1}$ of
$TM|_{V_3}$ by $E^u_{i1}=\im(id,\bar{\tau})$. Then
$E^u_{i1}|_{Q_3}=E^u_{i0}|_{Q_3}$,
$E^u_{i1}|_{W^u_{\delta'}(\Omega_i)}=E^u|_{W^u_{\delta'}(\Omega_i)}$,
and $Tf((E^u_{i1})_x)=(E^u_{i1})_{f(x)}$ if $x\in V_3\cap
f^{-1}(V_3)$.

Now consider the $C^\alpha$ and $d_f$-Lipschitz subbundle
$Tf^n(E^u_{i1})$ of $TM|_{f^n(V_3)}$, $n\in\mathbb{Z}$. If for
$n,m\in\mathbb{Z}$, $n<m$, $f^n(V_3)\cap f^m(V_3)\neq\emptyset$,
then for $x\in f^n(V_3)\cap f^m(V_3)$, $f^{-n}(x)\in V_3$,
$f^{-m}(x)\in V_3$. Since $V_3$ is unrevisited, $f^{-p}(x)\in V_3$
for $n\leq p\leq m$. So for $n+1\leq p\leq m$, $f^{-p}(x)\in V_3\cap
f^{-1}(V_3)$ and then
$Tf^p(E^u_{i1})_x=Tf^{p-1}(Tf((E^u_{i1})_{f^{-p}(x)}))
=Tf^{p-1}((E^u_{i1})_{f^{-p+1}(x)})=Tf^{p-1}(E^u_{i1})_x$. So
$Tf^n(E^u_{i1})_x=Tf^m(E^u_{i1})_x$, and then the bundles
$Tf^n(E^u_{i1})$ $(n\in\mathbb{Z})$ patch together to a
$Tf$-invariant subbundle $E^u_{i2}$ of $TM|_{\mathcal{O}(V_3)}$. It
is obviously continuous since $f^n(V_3)$ is open.

We have $T_xM=(E^u_{i2})_x+E^s_x$ for $x\in\mathcal{O}(V_3)$, as
this holds for $x\in W^u_{\delta'}(\Omega_i)\cup Q_3$, $E^u_{i2}$
and $E^s$ are $Tf$-invariant, and
$\mathcal{O}(V_3)=\mathcal{O}(W^u_{\delta'}(\Omega_i))\cup
\mathcal{O}(Q_3)$.

We prove that $(E^u_{i2})_x\subset(E^u_j)_x$ for every $j<i$ and
$x\in\mathcal{O}^-(V_3)\cap\mathcal{O}^+(U_j)$. Let
$x\in\mathcal{O}^-(V_3)\cap\mathcal{O}^+(U_j)$. Then $x\notin
W^u(\Omega_i)$ and then $x\in\mathcal{O}(Q_3)$. Since
$(E^u_{i2})_x\subset(E^u_j)_x$ for $x\in Q_3\cap\mathcal{O}^+(U_j)$,
it also holds for $x\in\mathcal{O}(Q_3)\cap\mathcal{O}^+(U_j)$ by
the $Tf$-invariance of $E^u_{i2}$ and $E^u_j$.

Finally, let $U_i$ be an open neighborhood of $\Omega_i$ such that
$\overline{U_i}\subset V_3$. Let $N>0$. We prove that $E^u_{i2}$ is
$C^\alpha$ and $d_f$-Lipschitz on $\bigcup_{n=-N}^{N}f^n(U_i)$.
Consider the Grassmanian bundle $\mathcal{G}$ over $M$ consisting of
all $\rank(E^u_{i2})$-dimensional subspaces of the tangent spaces of
$M$. Then $E^u_{i2}$ can be viewed as a $Tf$-invariant continuous
section $s$ of $\mathcal{G}|_{\mathcal{O}(V_3)}$ which is $C^\alpha$
and $d_f$-Lipschitz on each $f^n(V_3)$. Embed the compact manifold
$\mathcal{G}$ into some $\mathbb{R}^{N'}$. Then $s$ can be viewed as
a bounded map $s:\mathcal{O}(V_3)\rightarrow \mathbb{R}^{N'}$. Since
$\overline{f^n(U_i)}\subset f^n(V_3)$ for all $n$, there exists
$d_1>0$ such that for $-N\leq n\leq N$, $x\in f^n(U_i), y\in M,
d(x,y)<d_1$ imply that $y\in f^n(V_3)$. Let $K'>0$ be such that
$|s(x)-s(y)|\leq K'\rho_f(x,y)$ for $x,y\in f^n(V_3)$, $-N\leq n\leq
N$. We prove that $s$ is $C^\alpha$ and $d_f$-Lipschitz on
$\bigcup_{n=-N}^{N}f^n(U_i)$. Let
$x,y\in\bigcup_{n=-N}^{N}f^n(U_i)$. Since $s$ is bounded, we may
assume that $d(x,y)<d_1$. Suppose $x\in f^n(U_i)$. Then $y\in
f^n(V_3)$. Hence $|s(x)-s(y)|\leq K'\rho_f(x,y)$. So the
neighborhood $U_i$ of $\Omega_i$ and the bundle
$E^u_i=E^u_{i2}|_{\mathcal{O}(U_i)}$ satisfy the conditions
(i')--(iv'). The proof of Theorem \ref{T:1} is finished.
\end{proof}

\section{Existence of right inverses}

Let $f$ be a $C^2$ diffeomorphism of a compact manifold $M$,
$\alpha\in(0,1)$. Let $\mathfrak{X}^0(M)$ denote the Banach space of
continuous vector fields on $M$ with the $C^0$ norm $\|\cdot\|$, and
let $\mathfrak{X}^\alpha_f(M)$ be the subspace of
$\mathfrak{X}^0(M)$ consisting of $C^\alpha$ and $d_f$-Lipschitz
vector fields. As in the previous section, suppose $M$ is
isometrically embedded into some Euclidian space $\mathbb{R}^N$. For
$\eta\in\mathfrak{X}^\alpha_f(M)$, denote
$$L_\alpha(\eta)=\sup_{x,y\in M,x\neq y}\frac{|\eta(x)-\eta(y)|}{d(x,y)^\alpha},$$
$$L_f(\eta)=\sup_{x,y\in M,x\neq y}\frac{|\eta(x)-\eta(y)|}{d_f(x,y)}.$$ Then
$\mathfrak{X}^\alpha_f(M)$, being endowed with the norm
$$\|\eta\|_{\alpha,f}=\max\{\|\eta\|,L_\alpha(\eta),L_f(\eta)\},$$ is a
Banach space. For $\eta\in\mathfrak{X}^0(M)$, define the vector
field $f_\sharp(\eta)$ on $M$ by
$$f_\sharp(\eta)(x)=Tf(\eta(f^{-1}(x))).$$ Then $f_\sharp(\eta)$ is
in $\mathfrak{X}^0(M)$, and in $\mathfrak{X}^\alpha_f(M)$ if
$\eta\in\mathfrak{X}^\alpha_f(M)$.

The following theorem is motivated by \cite[Theorem B]{R1} and
\cite[Section 8]{R2}.

\begin{theorem}\label{T:2}
Let $f$ be a $C^2$ diffeomorphism of $M$ satisfying Axiom A and the
strong transversality condition, $\lambda, l$ be as in Theorem
\ref{T:1}. Suppose $\alpha\in(0,1)$ satisfies $\lambda l^\alpha<1$.
Then there exists a continuous linear operator $J$ on $\mathfrak{X}^0(M)$ such that\\
(i) $J$ is a right inverse of  $1-f_\sharp$;\\
(ii) $J$ maps $\mathfrak{X}^\alpha_f(M)$ into
$\mathfrak{X}^\alpha_f(M)$ and restricts to a continuous linear
operator on $\mathfrak{X}^\alpha_f(M)$ with respect to the norm
$\|\cdot\|_{\alpha,f}$.
\end{theorem}

\begin{proof}
Choose $\lambda<\lambda'<\rho=\kappa\lambda'<1$ such that $\rho
l^\alpha<1$. Let $\Omega=\Omega_1\cup\cdots\cup\Omega_k$ be the
spectral decomposition ordered as in Theorem \ref{T:1}. Let $U_i$ be
an open neighborhood of $\Omega_i$, $E^\sigma_i$ be two
$Tf$-invariant subbundles of $TM|_{\mathcal{O}(U_i)}$ satisfying
(i)--(iv) in Theorem \ref{T:1} for the above $\lambda'$,
$\sigma=u,s$. It is well known that
$\bigcup_{i=1}^k\mathcal{O}(U_i)=M$. So there exists $N>0$ such that
$\{\bigcup_{n=-N}^{N}f^n(U_1),\cdots,\bigcup_{n=-N}^{N}f^n(U_k)\}$
is a cover of $M$. Shrinking $U_i$ if necessary, we may assume that
they are unrevisited. Then it is easy to see that for every $x\in
M$, the set $\{n\in\mathbb{Z}:f^n(x)\notin\bigcup_{i=1}^kU_i\}$
contains at most $n_0=2kN$ elements. Let $\theta_1,\cdots,\theta_k$
be a smooth partition of unity subordinate to the above cover. For
$\eta\in\mathfrak{X}^0(M)$, let
$\eta_{i\sigma}=P_{E^\sigma_i}(\theta_i\eta)$, where $P_{E^s_i}$
(resp. $P_{E^u_i}$) is the projection of $TM_{\mathcal{O}(U_i)}$
onto $E^s_i$ (resp. $E^u_i$) along $E^u_i$ (resp. $E^s_i$), and
define $J_{is}(\eta)=\sum_{n=0}^{+\infty}f^n_\sharp(\eta_{is})$,
$J_{iu}(\eta)=-\sum_{n=1}^{+\infty}f^{-n}_\sharp(\eta_{iu})$,
$J(\eta)=\sum_{\sigma=u,s}\sum_{i=1}^kJ_{i\sigma}(\eta)$. Robbin
\cite{R1} proved that these series converge uniformly, and then $J$
is a continuous right inverse of $1-f_\sharp$. We prove in the
following that $J$ maps $\mathfrak{X}^\alpha_f(M)$ into
$\mathfrak{X}^\alpha_f(M)$ and restricts to a continuous linear
operator on $\mathfrak{X}^\alpha_f(M)$. As in \cite{R1}, it is
sufficient to prove this property for each $J_{is}$.

Let $\eta\in\mathfrak{X}^\alpha_f(M)$. Fix $i=1,\cdots,k$, and
denote $\zeta=\eta_{is}=P_{E^s_i}(\theta_i\eta)$. Then
$\supp(\zeta)\subset\bigcup_{n=-N}^{N}f^n(U_i)$ and
$\zeta(x)\in(E^s_i)_x$ for $x\in\bigcup_{n=-N}^{N}f^n(U_i)$. Since
$E^s_i$ and $E^u_i$ are $C^\alpha$ and $d_f$-Lipschitz on
$\bigcup_{n=-N}^{N}f^n(U_i)$, $\zeta\in\mathfrak{X}^\alpha_f(M)$.
Let $K=(\frac{\|Tf\|}{\rho})^{2n_0+N}$. It is proved in
\cite[Section 6]{R1} that
\begin{equation}\label{E:vecC0}
|f_\sharp^n(\zeta)(x)|\leq
(\frac{\|Tf\|}{\rho})^{n_0+N}\rho^n|\zeta(f^{-n}(x))|\leq
K\rho^n|\zeta(f^{-n}(x))|
\end{equation}
for all $x\in M$ and $n\geq0$ (note that we always have
$\|Tf\|>\rho$). Hence
\begin{equation}\label{E:fieldC0}
\|f_\sharp^n(\zeta)\|\leq K\rho^n\|\zeta\|
\end{equation}
for all $n\geq0$. Let
$$C=\|Tf\|\max\{L_\alpha(P_{E^s_j}|_{U_j}):1\leq j\leq
k\}+L_\alpha(Tf),$$ where $L_\alpha(Tf)$ is the H\"{o}lder constant
of $Tf$ as a map $x\mapsto T_xf$ for $x\in M$,
$L_\alpha(P_{E^s_j}|_{U_j})$ is the H\"{o}lder constant of
$P_{E^s_j}|_{U_j}$ as a map $x\mapsto P_{(E^s_j)_x}$ for $x\in U_j$.
We prove that
\begin{equation}\label{E:vecHol}
L_\alpha(f_\sharp^n(\zeta))\leq K(\rho l^\alpha)^nL_\alpha(\zeta)
+C'((\rho l^\alpha)^n-\rho^n)\|\zeta\|
\end{equation}
for all $n\geq0$, where $C'=\frac{CK^2l^\alpha}{\rho(l^\alpha-1)}$.

We first prove some inequalities on individual tangent vectors. Let
$p,q\in M$, $v_p\in T_pM, v_q\in T_qM$. Then
\begin{align}\label{E:v-bad}
&\notag|T_pf(v_p)-T_qf(v_q)|\\
\leq&|T_pf(v_p-v_q)|+|(T_pf-T_qf)(v_q)|\\
\leq&\notag\|Tf\| |v_p-v_q|+L_\alpha(Tf)|v_q|d(p,q)^\alpha\\
\leq&\notag\|Tf\| |v_p-v_q|+C|v_q|d(p,q)^\alpha.
\end{align}
Recall that a smooth adapted Riemannian metric on $M$ can be
obtained by approximating a $C^0$ adapted metric for which the
bundles $E^u$ and $E^s$ are mutually orthogonal on $\Omega$. So
after choosing a better approximation of the $C^0$ metric and
shrinking the $U_j$'s, we may assume that for each $j$,
$\|P_{(E^s_j)_p}\|\leq\kappa$ for every $p\in U_j$, where $\kappa>1$
is as in the beginning of the proof. So for $p,q,v_p,v_q$ as above,
if moreover we have $p,q\in U_j$ for some $j$, and $v_p\in(E^s_j)_p,
v_q\in(E^s_j)_q$, then
\begin{align}\label{E:v-good}
&\notag|T_pf(v_p)-T_qf(v_q)|\\
\leq&|T_pfP_{(E^s_j)_p}(v_p-v_q)|+|T_pf(P_{(E^s_j)_p}-P_{(E^s_j)_q})(v_q)|+|(T_pf-T_qf)(v_q)|\\
\leq&\notag\lambda'\kappa|v_p-v_q|+(\|Tf\|L_\alpha(P_{E^s_j}|_{U_j})+L_\alpha(Tf))|v_q|d(p,q)^\alpha\\
\leq&\notag\rho|v_p-v_q|+C|v_q|d(p,q)^\alpha.
\end{align}

Now we prove \eqref{E:vecHol}. Let $x,y\in M, n\geq0$. If one of
$f^{-n}(x)$ and $f^{-n}(y)$ does not belong to
$\bigcup_{n=-N}^{N}f^n(U_i)$, say
$f^{-n}(x)\notin\bigcup_{n=-N}^{N}f^n(U_i)$, then by
\eqref{E:vecC0}, we have
\begin{align*}
&|f_\sharp^n(\zeta)(x)-f_\sharp^n(\zeta)(y)|\\
=&|f_\sharp^n(\zeta)(y)|\\
\leq& K\rho^n|\zeta(f^{-n}(y))|\\
=&K\rho^n|\zeta(f^{-n}(x))-\zeta(f^{-n}(y))|\\
\leq&K(\rho l^\alpha)^nL_\alpha(\zeta)d(x,y)^\alpha.
\end{align*}
So \eqref{E:vecHol} holds in this case. Suppose
$f^{-n}(x),f^{-n}(y)\in\bigcup_{n=-N}^{N}f^n(U_i)$. Let $1\leq m\leq
n$. Then by letting $p=f^{-m}(x), q=f^{-m}(y),
v_p=f_\sharp^{n-m}(\zeta)(f^{-m}(x))$,
$v_q=f_\sharp^{n-m}(\zeta)(f^{-m}(y))$ in \eqref{E:v-bad} and using
\eqref{E:fieldC0}, we get
\begin{align}\label{E:f-bad}
&\notag|f_\sharp^{n-m+1}(\zeta)(f^{-m+1}(x))-f_\sharp^{n-m+1}(\zeta)(f^{-m+1}(y))|\\
=&\notag|T_{f^{-m}(x)}f(f_\sharp^{n-m}(\zeta)(f^{-m}(x)))-T_{f^{-m}(y)}f(f_\sharp^{n-m}(\zeta)(f^{-m}(y)))|\\
\leq&\|Tf\|
|f_\sharp^{n-m}(\zeta)(f^{-m}(x))-f_\sharp^{n-m}(\zeta)(f^{-m}(y))|\\
&\notag+C
\|f_\sharp^{n-m}(\zeta)\|d(f^{-m}(x),f^{-m}(y))^\alpha\\
\leq&\notag\|Tf\|
|f_\sharp^{n-m}(\zeta)(f^{-m}(x))-f_\sharp^{n-m}(\zeta)(f^{-m}(y))|\\
&\notag+CK\rho^{n-m}l^{m\alpha}\|\zeta\|d(x,y)^\alpha.
\end{align}
If moreover $f^{-m}(x), f^{-m}(y)\in U_j$ for some $j\geq i$, then
$f_\sharp^{n-m}(\zeta)(f^{-m}(x))\in(E^s_j)_{f^{-m}(x)}$,
$f_\sharp^{n-m}(\zeta)(f^{-m}(y))\in(E^s_j)_{f^{-m}(y)}$. By
\eqref{E:v-good} and \eqref{E:fieldC0},
\begin{align}\label{E:f-good}
&\notag|f_\sharp^{n-m+1}(\zeta)(f^{-m+1}(x))-f_\sharp^{n-m+1}(\zeta)(f^{-m+1}(y))|\\
\leq&\notag\rho
|f_\sharp^{n-m}(\zeta)(f^{-m}(x))-f_\sharp^{n-m}(\zeta)(f^{-m}(y))|\\
&+C
\|f_\sharp^{n-m}(\zeta)\|d(f^{-m}(x),f^{-m}(y))^\alpha\\
\leq&\notag\rho
|f_\sharp^{n-m}(\zeta)(f^{-m}(x))-f_\sharp^{n-m}(\zeta)(f^{-m}(y))|\\
&\notag+CK\rho^{n-m}l^{m\alpha}\|\zeta\|d(x,y)^\alpha.
\end{align}
For $1\leq m\leq n$, denote
$$\nu_m=
\begin{cases}
\rho, &\text{if $f^{-m}(x), f^{-m}(y)\in U_j$ for some $j\geq i$;}\\
\|Tf\|, &\text{otherwise.}\\
\end{cases}$$
Then by \eqref{E:f-bad} and \eqref{E:f-good}, we have
\begin{align}\label{E:f}
&\notag|f_\sharp^{n-m+1}(\zeta)(f^{-m+1}(x))-f_\sharp^{n-m+1}(\zeta)(f^{-m+1}(y))|\\
\leq&\nu_m
|f_\sharp^{n-m}(\zeta)(f^{-m}(x))-f_\sharp^{n-m}(\zeta)(f^{-m}(y))|\\
&\notag +CK\rho^{n-m}l^{m\alpha}\|\zeta\|d(x,y)^\alpha.
\end{align}
Since we have supposed that
$f^{-n}(x),f^{-n}(y)\in\bigcup_{n=-N}^{N}f^n(U_i)$, we have
$f^{-(n-N)}(x)$, $f^{-(n-N)}(y)\in\mathcal{O}^+(U_i)$. But
$\mathcal{O}^+(U_i)\cap \mathcal{O}^-(U_j)=\emptyset$ for $j<i$. So
each of the sets $\{1\leq m\leq n-N:f^{-m}(x)\notin
\bigcup_{j=i}^kU_j\}$ and $\{1\leq m\leq n-N:f^{-m}(y)\notin
\bigcup_{j=i}^kU_j\}$ consists of at most $n_0$ elements. Then for
all but at most $2n_0+N$ integers $m$ in $\{1,\cdots,n\}$,
$f^{-m}(x), f^{-m}(y)\in U_j$ for some $j\geq i$, that is, at most
$2n_0+N$ numbers $\nu_m (1\leq m\leq n)$ equal to $\|Tf\|$. So we
have $\nu_1\nu_2\cdots\nu_m\leq
(\frac{\|Tf\|}{\rho})^{2n_0+N}\rho^m=K\rho^m$. Then by \eqref{E:f},
we get
\begin{align*}
&|f_\sharp^{n}(\zeta)(x)-f_\sharp^{n}(\zeta)(y)|\\
\leq&\nu_1\nu_2\cdots\nu_n
|\zeta(f^{-n}(x))-\zeta(f^{-n}(y))|\\
&+CK(\rho^{n-1}l^\alpha+\nu_1\rho^{n-2}l^{2\alpha}+\nu_1\nu_2\rho^{n-3}l^{3\alpha}\\
&+\cdots+\nu_1\nu_2\cdots\nu_{n-1}l^{n\alpha})\|\zeta\|d(x,y)^\alpha\\
\leq& K\rho^n|\zeta(f^{-n}(x))-\zeta(f^{-n}(y))|\\
&+CK^2\rho^{n-1}(l^\alpha+l^{2\alpha}+\cdots+l^{n\alpha})\|\zeta\|d(x,y)^\alpha\\
\leq& K(\rho l^\alpha)^nL_\alpha(\zeta)
d(x,y)^\alpha+\frac{CK^2l^\alpha}{\rho(l^\alpha-1)}((\rho
l^\alpha)^n-\rho^n)\|\zeta\|d(x,y)^\alpha.
\end{align*}
This proves \eqref{E:vecHol}.

Recall that $\zeta=\eta_{is}$. By \eqref{E:vecHol}, we have
\begin{equation}\label{E:JHol}
L_\alpha(J_{is}(\eta))
\leq\sum_{n=0}^{+\infty}L_\alpha(f_\sharp^{n}(\eta_{is}))
\leq\frac{K}{1-\rho l^\alpha}L_\alpha(\eta_{is})+(\frac{C'}{1-\rho
l^\alpha}-\frac{C'}{1-\rho})\|\eta_{is}\|.
\end{equation}
Similarly, we can prove that
\begin{equation}\label{E:JLip}
L_f(J_{is}(\eta)) \leq AL_f(\eta_{is})+B\|\eta_{is}\|
\end{equation}
for some constant $A,B>0$ (see \cite[Section 6]{R1}). Since the
bundles $E^s_i$ and $E^u_i$ are $C^\alpha$ and $d_f$-Lipschitz on
$\bigcup_{n=-N}^{N}f^n(U_i)$, the operator on
$\mathfrak{X}^\alpha_f(M)$ which maps $\eta$ to $\eta_{is}$ is
continuous. So by \eqref{E:fieldC0}, \eqref{E:JHol} and
\eqref{E:JLip}, the operator $J_{is}$ maps
$\mathfrak{X}^\alpha_f(M)$ into $\mathfrak{X}^\alpha_f(M)$ and is
continuous on $\mathfrak{X}^\alpha_f(M)$. This proves the theorem.
\end{proof}

\section{Proof of Theorem \ref{T:main2}}

In this section we prove Theorem \ref{T:main2}. As indicated in the
introduction, Theorem \ref{T:main} follows from Theorem
\ref{T:main2}.

We first extract some analytical arguments in \cite{R1,R2} to the
following lemma, which can be viewed as a generalization of the
usual Implicit Function Theorem for Banach spaces.

\begin{lemma}\label{L:implicit}
Let $(X,\|\cdot\|)$ be a Banach space, $X'$ be a linear subspace of
$X$ with a complete norm $\|\cdot\|'$ such that the inclusion
$(X',\|\cdot\|')\hookrightarrow(X,\|\cdot\|)$ is continuous, and
such that the closed unit ball $\{x\in X':\|x\|'\leq1\}$ in $X'$ is
a closed subset of $(X,\|\cdot\|)$. Let $\mathcal{M}$ be a Banach
manifold, $f\in\mathcal{M}$, $\mathcal{U}$ be an open set in $X$
containing $0\in X$. Let $\Psi:\mathcal{M}\times
\mathcal{U}\rightarrow X$ be a $C^1$ map satisfying $\Psi(f,0)=0$
and $\Psi(\mathcal{M}\times(\mathcal{U}\cap X'))\subset X'$. Denote
by $A=D_2\Psi(f,0):X\rightarrow X$ the partial derivative of $\Psi$
at the point $(f,0)$ along the second
variable. Suppose\\
(1) $A(X')\subset X'$;\\
(2) $1-A$ has a continuous linear right inverse $J$ which maps $X'$
into $X'$ and restricts to a continuous linear operator on $X'$;\\
(3) for any $\varepsilon>0$, there exist a neighborhood
$\mathcal{M}_\varepsilon$ of $f$ in $\mathcal{M}$ and a neighborhood
$\mathcal{U}_\varepsilon$ of $0$ in $\mathcal{U}$ such that
$$\|\Psi(g,x)-A(x)\|'\leq\varepsilon(1+\|x\|')$$
for all $g\in\mathcal{M}_\varepsilon$, $x\in\mathcal{U}_\varepsilon\cap X'$.\\
Then for any neighborhood $\mathcal{V}\subset X'$ of $0$ in
$(X',\|\cdot\|')$, there exist a neighborhood $\mathcal{N}$ of $f$
in $\mathcal{M}$ and a map $c:\mathcal{N}\rightarrow\mathcal{V}$
such that\\
(i) $c(f)=0$;\\
(ii) $\Psi(g,c(g))=c(g)$ for all $g\in\mathcal{N}$;\\
(iii) as a map $\mathcal{N}\rightarrow X$, $c$ is $C^1$.
\end{lemma}

\begin{proof}
Denote the norm of $J$ as a operator on $X$ by $\|J\|$, and the norm
of $J|_{X'}$ as a operator on $X'$ by $\|J\|'$. Choose
$0<\varepsilon\leq1$ such that the closed ball
$B'(\varepsilon)=\{x\in X':\|x\|'\leq \varepsilon\}$ lies in
$\mathcal{V}$. By the condition (3) and the continuous
differentiability of $\Psi$, we may choose an open neighborhood
$\mathcal{N}$ of $f$ in $\mathcal{M}$ and $r>0$ such that the closed
ball $B(r)=\{x\in X:\|x\|\leq r\}$ lies in $\mathcal{U}$, and such
that
\begin{equation}\label{E:l1}
\|D_2\Psi(g,x)-A\|\leq\frac{1}{2\|J\|}
\end{equation}
for all
$g\in\mathcal{N}, x\in B(r)$, and
\begin{equation}\label{E:l2}
\|\Psi(g,x)-A(x)\|'\leq\frac{\varepsilon}{2\|J\|'}(1+\|x\|')
\end{equation}
for all
$g\in\mathcal{N}, x\in B(r)\cap X'$. By making $\mathcal{N}$
smaller, we may also assume that
\begin{equation}\label{E:l3}
\|\Psi(g,0)\|\leq\frac{r}{2\|J\|}
\end{equation}
for all $g\in\mathcal{N}$.

For $g\in\mathcal{N}$, define a map $R_g:B(r)\rightarrow X$ by
$$R_g(x)=J(\Psi(g,x)-A(x)).$$ Then for $x\in B(r)$, by \eqref{E:l1}, \eqref{E:l3}
and the Mean Value Theorem, we have
\begin{align*}
&\|R_g(x)\|\\
\leq&\|J\|(\|\Psi(g,0)\|+\|(\Psi(g,x)-A(x))-(\Psi(g,0)-A(0))\|)\\
\leq&\|J\|(\frac{r}{2\|J\|}+\frac{1}{2\|J\|}\|x\|)\\
\leq&r.
\end{align*}
So $R_g$ maps $B(r)$ into $B(r)$. For $x,y\in B(r)$, also by
\eqref{E:l1} and the Mean Value Theorem, we have
\begin{align*}
&\|R_g(x)-R_g(y)\|\\
\leq&\|J\|\|(\Psi(g,x)-A(x))-(\Psi(g,y)-A(y))\|\\
\leq&\|J\|\frac{1}{2\|J\|}\|x-y\|\\
=&\frac{1}{2}\|x-y\|.
\end{align*}
So $R_g$ is a contraction on $B(r)$. By the Contraction Principle,
there is a unique fixed point $c(g)$ of $R_g$ in $B(r)$. This means
that $(1-A)(c(g))=(1-A)(R_g(c(g)))=\Psi(g,c(g))-A(c(g))$. So
$\Psi(g,c(g))=c(g)$. It is obvious that $c(f)=0$.

We prove that $c(g)\in\mathcal{V}$. Let $x_n=R_g^n(0)\in B(r),
n\geq0$. Then $\|x_n-c(g)\|\rightarrow0$, and it is obvious by
induction that $x_n\in X'$. We have
$x_{n+1}=R_g(x_n)=J(\Psi(g,x_n)-A(x_n)$. By \eqref{E:l2}, we get
$\|x_{n+1}\|'\leq\frac{\varepsilon}{2}(1+\|x\|')$, which is
equivalent to
$$\|x_{n+1}\|'-\frac{\varepsilon}{2-\varepsilon}\leq
\frac{\varepsilon}{2}(\|x_{n}\|'-\frac{\varepsilon}{2-\varepsilon}).$$
By induction we easily get
$\|x_{n}\|'-\frac{\varepsilon}{2-\varepsilon}\leq0$ for all
$n\geq0$. Hence
$\|x_{n}\|'\leq\frac{\varepsilon}{2-\varepsilon}\leq\varepsilon$.
But the closed ball $B'(\varepsilon)$ in $X'$ is closed in $X$ and
$x_n\rightarrow c(g)$ in $X$. So $c(g)\in
B'(\varepsilon)\subset\mathcal{V}$.

The proof of the fact that $c$ as a map $\mathcal{N}\rightarrow X$
is $C^1$ is the same as the proof of the corresponding result in the
usual Implicit Function Theorem. We omit the details here.
\end{proof}

\begin{proof}[Proof of Theorem \ref{T:main2}]
The map $\Psi:\Diff\times C^0(M,M)\rightarrow C^0(M,M)$ between
Banach manifolds defined by $$\Psi(g,h)=ghf^{-1}$$ is $C^1$ (see,
for example, \cite{Fr}). Let $(\mathcal{U}_0,\varphi)$ be a
coordinate chart around the identity map $id$ in $C^0(M,M)$, where
the coordinate $\varphi:\mathcal{U}_0\rightarrow\mathfrak{X}^0(M)$
is provided by the exponential map associated with some Riemannian
metric on $M$, that is, $\varphi(h)(x)=\exp_x^{-1}(h(x))$. $\varphi$
maps the set of $C^\alpha$ and $d_f$-Lipschitz maps in
$\mathcal{U}_0$ onto
$\varphi(\mathcal{U}_0)\cap\mathfrak{X}^\alpha_f(M)$. Let
$\mathcal{U}\subset \mathcal{U}_0$ be an open neighborhood of $id$
in $C^0(M,M)$, $\mathcal{M}$ be an open neighborhood of $f$ in
$\Diff$, such that $\Psi(\mathcal{M}\times \mathcal{U})\subset
\mathcal{U}_0$. By abuse of language, we identify $\mathcal{U}_0$
with $\varphi(\mathcal{U}_0)$ via the coordinate $\varphi$. But we
denote an element in $\mathcal{U}_0$ by $h$ when we view it as a
map, and by $\eta$ if it is regarded as a vector field.

The partial derivative
$D_2\Psi(f,id):\mathfrak{X}^0(M)\rightarrow\mathfrak{X}^0(M)$ of
$\Psi$ at the point $(f,id)$ equals to $f_\sharp$. Since $f$ is
$C^2$, $D_2\Psi(f,id)$ maps $\mathfrak{X}^\alpha_f(M)$ into
$\mathfrak{X}^\alpha_f(M)$. By Theorem \ref{T:2}, $1-D_2\Psi(f,id)$
has a right inverse $J$ which restricts to a continuous linear
operator on $\mathfrak{X}^\alpha_f(M)$.

To apply Lemma \ref{L:implicit}, we need to verify the following two
conditions.\\
(1) The closed unit ball in $\mathfrak{X}^\alpha_f(M)$ is a closed
subset in $\mathfrak{X}^0(M)$;\\
(2) For every $\varepsilon>0$, there exist a neighborhood
$\mathcal{M}_\varepsilon\subset \mathcal{M}$ of $f$ and $\delta>0$
such that
$\|\Psi(g,\eta)-f_\sharp(\eta)\|_{\alpha,f}\leq\varepsilon(1+\|\eta\|_{\alpha,f})$
for all $g\in\mathcal{M}_\varepsilon, \eta\in
\mathcal{U}\cap\mathfrak{X}^\alpha_f(M)$ with $\|\eta\|<\delta$.

To prove (1), let $(\eta_n)_{n=1}^{\infty}$ be a sequence in the
closed unit ball in $\mathfrak{X}^\alpha_f(M)$, that is,
$\|\eta_n\|_{\alpha,f}=\max\{\|\eta_n\|,L_\alpha(\eta_n),L_f(\eta_n)\}\leq1$
for all $n$. Suppose $\eta\in\mathfrak{X}^0(M)$ such that
$\|\eta_n-\eta\|\rightarrow0$. Then $\|\eta\|\leq1$. By letting
$n\rightarrow\infty$ in the inequality
$\frac{|\eta_n(x)-\eta_n(y)|}{d(x,y)^\alpha}\leq1$, we get
$L_\alpha(\eta)\leq1$. Similarly, $L_f(\eta)\leq1$. So
$\|\eta\|_{\alpha,f}\leq1$. (1) is proved.

Denote $Q(g,\eta)=\Psi(g,\eta)-f_\sharp(\eta)$. Let
$\varepsilon'>0$. Then
\begin{equation}\label{E:b1}
\|Q(g,\eta)\|\leq\varepsilon'
\end{equation}
for $g$ sufficiently $C^1$ close to $f$ and $\|\eta\|$ sufficiently
small. By considering the partial differentials of the $C^1$ map
$\mathcal{M}\times TM\rightarrow TM,
(g,x,v)\mapsto(f(x),\exp_{f(x)}^{-1}(g(\exp_x(v))))$ along the
directions of $x$ and $v$ (see \cite[Lemma 3.2, Lemma 3.4]{R1} or
\cite[Lemma 8.4]{R2}), we have
$$|Q(g,\eta)(x)-Q(g,\eta)(y)|\leq\varepsilon'(d(f^{-1}(x),f^{-1}(y))+|\eta(f^{-1}(x))-\eta(f^{-1}(y))|)$$
whenever $g$ is sufficiently $C^1$ close to $f$ and $\|\eta\|$ is
sufficiently small, from which we easily get
\begin{equation}\label{E:b2}
L_\alpha(Q(g,\eta))\leq\varepsilon'(l\diam(M)^{1-\alpha}+l^\alpha
L_\alpha(\eta)),
\end{equation}
\begin{equation}\label{E:b3}
L_f(Q(g,\eta))\leq\varepsilon'(1+L_f(\eta))
\end{equation}
for such $g$ and $\eta$ if $\eta\in\mathfrak{X}^\alpha_f(M)$, where
$\diam(M)$ is the diameter of $M$. By \eqref{E:b1}, \eqref{E:b2} and
\eqref{E:b3}, we have
$$\|Q(g,\eta)\|_{\alpha,f}\leq\varepsilon'\max\{l\diam(M)^{1-\alpha},l^\alpha\}(1+\|\eta\|_{\alpha,f})$$
for $g$ sufficiently $C^1$ close to $f$ and
$\eta\in\mathfrak{X}^\alpha_f(M)$ with $\|\eta\|$ sufficiently
small. This proves (2).

Let $\mathcal{V}$ be a $C^\alpha$ neighborhood of $id$ in
$C^\alpha(M,M)$ as in Theorem \ref{T:main2}. Then we may choose a
neighborhood $\mathcal{V}_f$ of $id$ in the Banach manifold
$C^\alpha_f(M,M)$ of $C^\alpha$ and $d_f$-Lipschitz maps on $M$ such
that $\mathcal{V}_f\subset\mathcal{V}$, and such that elements in
$\mathcal{V}_f$ are sufficiently $d_f$-Lipschitz close to the
identity. Applying Lemma \ref{L:implicit} to the map $\Psi$, we get
a $C^1$ neighborhood $\mathcal{N}$ of $f$ in
$\mathcal{M}\subset\Diff$ and a function $c:\mathcal{N}\rightarrow
\mathcal{V}_f$ with $c(f)=id$ such that
$\Psi(g,c(g))=gc(g)f^{-1}=c(g)$ for every $g\in\mathcal{N}$, and $c$
is $C^1$ as a map $\mathcal{N}\rightarrow C^0(M,M)$. It is easy to
show that if $c(g)$ is sufficiently $d_f$-Lipschitz close to the
identity, then $c(g)$ is a homeomorphism (see \cite{R1,R2}). So
$g=c(g)fc(g)^{-1}$. This proves Theorem \ref{T:main2}.
\end{proof}

\end{document}